# GLOBAL EXPONENTIAL OBSERVERS
# FOR TWO CLASSES OF NONLINEAR SYSTEMS

**Iasson Karafyllis[*] and Costas Kravaris[**]**

[*]**Department of Environmental Engineering, Technical University of Crete,
73100, Chania, Greece (email: ikarafyl@enveng.tuc.gr )**

[**]**Department of Chemical Engineering, University of Patras
26500 Patras, Greece (email: kravaris@chemeng.upatras.gr )**

**Abstract**
This paper develops sufficient conditions for the existence of global exponential observers for two classes of nonlinear systems: (i) the class of systems with a globally asymptotically stable compact set, and (ii) the class of systems that evolve on an open set. In the first class, the derived continuous-time observer also leads to the construction of a robust global sampled-data exponential observer, under additional conditions. Two illustrative examples of applications of the general results are presented, one is a system with monotone nonlinearities and the other is the chemostat system.

## 1. Introduction

One of the biggest challenges of mathematical control theory has been the problem of constructing state observers for nonlinear systems. This problem has attracted a lot of attention in the literature in the past decades; it has been approached with a variety of methods and from a variety of points of view (see for instance [1,2,4,5,7,8,11,16,17,18] and references therein). In this work, we focus on nonlinear forward complete systems of the form:

$$\dot{x} = f(x,u), x \in \Re^n, u \in U \qquad (1.1)$$

where $U \subseteq \Re^m$ is a non-empty set, $f : \Re^n \to \Re^n$ is a smooth vector field and the output is given by

$$y = h(x) \qquad (1.2)$$

where $h : \Re^n \to \Re^k$ is a smooth mapping. The aim is to construct global exponential observers.

Available methods for global exponential observers include high-gain observers for globally Lipschitz systems ([8]) as well as circle-criterion observers, primarily for nonlinear systems with monotone nonlinearities ([1,7]). In transformation-based observers, originally developed in local form in [11] and subsequently in global form in [2], the system is mapped to a linear system, and the design of the observer is performed in transformed coordinates, where exponential convergence is imposed.

In this work, we present sufficient conditions for the existence of exponential observers for two important classes of nonlinear systems, which are not covered by the above methods:

1) Nonlinear systems with an asymptotically stable compact set

2) Nonlinear systems evolving on open sets

For both classes of systems, the proposed construction of the global exponential observer starts with a "candidate observer", which is subsequently modified by adding a correction term, in order to satisfy appropriate Lyapunov inequalities. It should be emphasized that explicit formulae for the observers are provided in each case and therefore the control practitioner can directly apply the results of the paper.

In Section 2, where we study the first class of systems, the "candidate observer" is a local observer over a certain compact set, whereas the correction term forces the trajectory to enter the compact set in finite time. The derived continuous-time observer can also lead to the construction of a robust global sampled-data exponential observer, under additional conditions. The sampled-data exponential observer is robust with respect to perturbations of the sampling schedule and with respect to measurement errors (see also [9,13,14] for sampled-data observers).

Section 3 studies the second class of systems, with the property of evolving on an open proper subset of $\Re^n$. Here, the "candidate observer" does not guarantee that the observer trajectories lie within the open set, and this is accomplished by adding an appropriate correction term. The design of the correction term is performed after transforming the system through an appropriate smooth injective map that maps the open set onto $\Re^n$, even though exponential convergence is enforced in the original coordinates. The results of Section 3 are important because for many classes of systems the state evolves in an open set (for example, biological systems usually evolve in the open first quadrant). However, there is another reason that motivates the results of Section 3. If a change of coordinates $X = \Phi(x)$ can be found, where $\Phi : \Re^n \to \Re^n$ is a smooth injective mapping satisfying $D\Phi(x)f(x,u) = A\Phi(x) + b(h(x),u)$ for all $x \in \Re^n$ for certain Hurwitz matrix $A \in \Re^{n \times n}$ and certain mapping $b : h(\Re^n) \times U \to \Re^n$, then the mapping $\Phi : \Re^n \to \Re^n$ can be used for the design of an observer for (1.1), (1.2) under additional hypotheses (see [2,11]). The results of Section 3 show that we do not have to assume that $\Phi(\Re^n) = \Re^n$ (i.e., $\Phi : \Re^n \to \Re^n$ is onto): instead we can require that $A = \Phi(\Re^n)$ is an open set and apply Theorem 3.1.

Finally, in Section 4, we present two illustrative examples of application of the general results. The first example is a system with monotone nonlinearities, and we apply the results of Section 2 to derive a global exponential observer, first under continuous-time measurements and subsequently under sampled measurements. The second example is a bioreactor, following the chemostat model, with positive state variables evolving on the open first quadrant of $\Re^2$. Applying the results of Section 3 leads to a global exponential observer, with positive state estimates.

**Notation.** Throughout this paper, we adopt the following notation:

* $\Re_+ := [0,+\infty)$.
* By $C^0(A;\Omega)$, we denote the class of continuous functions on $A \subseteq \Re^n$, which take values in $\Omega \subseteq \Re^m$. By $C^k(A;\Omega)$, where $k \geq 1$ is an integer, we denote the class of functions on $A \subseteq \Re^n$ with continuous derivatives of order $k$, which take values in $\Omega \subseteq \Re^m$.
* By $\text{int}(A)$, we denote the interior of the set $A \subseteq \Re^n$.
* For a vector $x \in \Re^n$, we denote by $x'$ its transpose and by $|x|$ its Euclidean norm. $A' \in \Re^{n \times m}$ denotes the transpose of the matrix $A \in \Re^{m \times n}$ and $|A|$ denotes the induced norm of the matrix $A \in \Re^{m \times n}$, i.e., $|A| = \sup\{|Ax| : x \in \Re^m, |x| = 1\}$.
* A function $V : \Re^n \to \Re_+$ will be called positive definite if $V(0) = 0$ and $V(x) > 0$ for all $x \neq 0$. A function $V : \Re^n \to \Re_+$ will be called radially unbounded if the sets $\{x \in \Re^n : V(x) \leq M\}$ are either empty or bounded for all $M \geq 0$.
* For a function $V \in C^1(A;\Re)$, the gradient of $V$ at $x \in A \subseteq \Re^n$, denoted by $\nabla V(x)$, is the row vector
$$\nabla V(x) = \left[ \frac{\partial V}{\partial x_1}(x) \quad \ldots \quad \frac{\partial V}{\partial x_n}(x) \right].$$

## 2. Systems with a Globally Asymptotically Stable Compact Set

Consider the forward complete system (1.1), (1.2). Our main hypothesis in this section guarantees that there exists a compact set which is robustly globally asymptotically stable (the adjective *robust* means uniformity to all measurable and locally essentially bounded inputs $u : \Re_+ \to U$ ).

**(H1)** *There exist a radially unbounded (but not necessarily positive definite) function $V \in C^2(\Re^n; \Re_+)$, a positive definite function $W \in C^1(\Re^n; \Re_+)$ and a constant $R > 0$ such that the following inequality holds for all $(x, u) \in \Re^n \times U$ with $V(x) \geq R$*

$$\nabla V(x) f(x, u) \leq -W(x) \tag{2.1}$$

Indeed, hypothesis (H1) guarantees that for every initial condition $x(0) \in \Re^n$ and for every measurable and locally essentially bounded input $u : \Re_+ \to U$ the solution $x(t)$ of (1.1) enters the compact set $S = \{x \in \Re^n : V(x) \leq R\}$ after a finite transient period, i.e., there exists $T \in C^0(\Re^n; \Re_+)$ such that $x(t) \in S$, for all $t \geq T(x(0))$. Moreover, notice that the compact set $S = \{x \in \Re^n : V(x) \leq R\}$ is positively invariant. This fact is guaranteed by the following lemma, which is proved at the Appendix.

**Lemma 2.1:** *Consider system (1.1) under hypothesis (H1). Then there exists $T \in C^0(\Re^n; \Re_+)$ such that for every $x_0 \in \Re^n$ and for every measurable and locally essentially bounded input $u : \Re_+ \to U$ the solution $x(t) \in \Re^n$ of (1.1) with initial condition $x(0) = x_0$ and corresponding to input $u : \Re_+ \to U$ satisfies $V(x(t)) \leq \max(V(x_0), R)$ for all $t \geq 0$ and $V(x(t)) \leq R$ for all $t \geq T(x_0)$.*

Our second hypothesis guarantees that we are in a position to construct an appropriate local exponential observer for system (1.1), (1.2).

**(H2)** *There exist a symmetric and positive definite matrix $P \in \Re^{n \times n}$, constants $\mu > 0$, $b > R$ and a smooth mapping $k : \Re^n \times h(\Re^n) \times U \to \Re^n$ with $k(\xi, y, u) = 0$ for all $(\xi, y, u) \in \Re^n \times h(\Re^n) \times U$ with $h(\xi) = y$ such that the following inequality holds*

$$(\xi - x)' P \big( f(\xi, u) + k(\xi, h(x), u) - f(x, u) \big) \leq -\mu |\xi - x|^2, \text{ for all } u \in U, \xi, x \in \Re^n \text{ with } V(\xi) \leq b \text{ and } V(x) \leq R \tag{2.2}$$

Indeed, hypothesis (H2) in conjunction with hypothesis (H1) guarantees that for every $x(0) \in S = \{x \in \Re^n : V(x) \leq R\}$ and for every measurable and locally essentially bounded input $u : \Re_+ \to U$ the solution of system (1.1), (1.2) with

$$\dot{\xi} = f(\xi, u) + k(\xi, y, u) \tag{2.3}$$

will satisfy an estimate of the form $|\xi(t) - x(t)| \leq M \exp(-\sigma t)|\xi(0) - x(0)|$, for all $t \geq 0$ for appropriate constants $M, \sigma > 0$, provided that the initial estimation error $|\xi(0) - x(0)|$ is sufficiently small. This is why system (2.3) is termed as "a local exponential observer". The reader should notice that hypothesis (H2) holds automatically for nonlinear systems of the form

$$\begin{aligned}
\dot{x}_1 &= f_1(x_1) + x_2 \\
\dot{x}_2 &= f_2(x_1, x_2) + x_3 \\
&\vdots \\
\dot{x}_n &= f_n(x_1, \ldots, x_n) + u \\
y &= x_1
\end{aligned} \tag{2.4}$$

for every $b > R > 0$ and for every non-empty set $U \subseteq \Re^m$, where $f_i : \Re^i \to \Re$ ($i = 1, \ldots, n$) are smooth mappings.

In order to be able to construct a nonlinear exponential observer for system (1.1), (1.2) we need an additional technical hypothesis.

**(H3)** *There exist constants* $c \in (0,1)$, $R \leq a < b$ *such that the following inequality holds:*

$$\nabla V(\xi)(f(\xi,u)+k(\xi,h(x),u)) \leq -W(\xi) + (1-c)|\nabla V(\xi)|^2 \frac{(\xi-x)'P(f(\xi,u)+k(\xi,h(x),u)-f(x,u))}{\nabla V(\xi)P(\xi-x)}$$

for all $u \in U$, $\xi,x \in \Re^n$ with $a < V(\xi) \leq b$, $\nabla V(\xi)P(\xi-x) < 0$ and $V(x) \leq R$ \hfill (2.5)

Hypothesis (H3) imposes constraints for the evolution of the trajectories of the local observer (2.3). Indeed, inequality (2.5) imposes a bound on the derivative of the Lyapunov function $V \in C^1(\Re^n;\Re_+)$ along the trajectories of the local observer (2.3) for specific regions of the state space.

We are now ready to state and prove the main result of the present section.

**Theorem 2.2:** *Consider system (1.1), (1.2) under hypotheses (H1-3). Define the locally Lipschitz mapping* $\hat{k}: \Re^n \times h(\Re^n) \times U \to \Re^n$:

$$\hat{k}(\xi,y,u) := k(\xi,y,u), \text{ for all } (\xi,y,u) \in \Re^n \times h(\Re^n) \times U \text{ with } V(\xi) \leq R \tag{2.6}$$

$$\hat{k}(\xi,y,u) := k(\xi,y,u) - \frac{\varphi(\xi,y,u)}{|\nabla V(\xi)|^2}(\nabla V(\xi))', \text{ for all } (\xi,y,u) \in \Re^n \times h(\Re^n) \times U \text{ with } V(\xi) > R \tag{2.7}$$

*where* $\varphi: \Re^n \times h(\Re^n) \times U \to \Re_+$ *is defined by*

$$\varphi(\xi,y,u) := \max\left(0, \nabla V(\xi)f(\xi,u)+W(\xi)+p(V(\xi))\nabla V(\xi)k(\xi,y,u)\right) \tag{2.8}$$

*and* $p: \Re_+ \to [0,1]$ *is an arbitrary locally Lipschitz function that satisfies* $p(s) = 1$ *for all* $s \geq b$ *and* $p(s) = 0$ *for all* $s \leq a$. *Then there exists* $M \in C^0(\Re^n \times \Re^n;\Re_+)$ *and* $\sigma > 0$ *such that for every measurable and locally essentially bounded input* $u: \Re_+ \to U$ *and for every* $(x_0,\xi_0) \in \Re^n \times \Re^n$ *the solution* $(x(t),\xi(t)) \in \Re^n \times \Re^n$ *of (1.1), (1.2) with*

$$\dot{\xi}(t) = f(\xi,u) + \hat{k}(\xi,y,u) \tag{2.9}$$

*initial condition* $(x(0),\xi(0)) = (x_0,\xi_0)$ *corresponding to input* $u: \Re_+ \to U$ *satisfies*

$$|\xi(t)-x(t)| \leq M(x_0,\xi_0)\exp(-\sigma t), \text{ for all } t \geq 0 \tag{2.10}$$

**Remark 2.3:**
**(a)** Theorem 2.2 shows that under hypotheses (H1-3), a "correction term" is needed in order to be able to construct a global exponential observer for system (1.1), (1.2). The "correction term" $-\frac{\varphi(\xi,y,u)}{|\nabla V(\xi)|^2}(\nabla V(\xi))'$ becomes active in the region $V(\xi) > a$ and its main task is to guarantee the validity of the differential inequality $\nabla V(\xi)(f(\xi,u)+\hat{k}(\xi,h(x),u)) \leq -W(\xi)$ for all $(\xi,x,u) \in \Re^n \times \Re^n \times U$ with $V(\xi) \geq b$. The previous differential inequality in conjunction with Lemma 2.1 guarantees that the solution enters an appropriate compact set in finite time and in this appropriate compact set the local exponential observer works.

**(b)** Inequality (2.2) guarantees that hypothesis (H3) holds provided that there exist constants $c \in (0,1)$, $R \leq a < b$ such that the following inequality holds for all $u \in U$, $\xi,x \in \Re^n$ with $a < V(\xi) \leq b$ and $V(x) \leq R$:

$$\nabla V(\xi)(f(\xi,u)+k(\xi,h(x),u)) \leq -W(\xi) + \frac{\mu(1-c)}{|P|}|\nabla V(\xi)||\xi-x| \tag{2.11}$$

**Proof of Theorem 2.2:** First notice that for all $(\xi, x, u) \in \Re^n \times \Re^n \times U$ with $V(\xi) \geq b$ the following inequality holds:

$$\nabla V(\xi)(f(\xi, u) + \hat{k}(\xi, h(x), u)) \leq -W(\xi) \tag{2.12}$$

Indeed, definition (2.7) implies $\nabla V(\xi)(f(\xi,u) + \hat{k}(\xi, h(x), u)) = \nabla V(\xi)(f(\xi,u) + k(\xi, h(x), u)) - \varphi(\xi, h(x), u)$. By distinguishing the cases $\nabla V(\xi) f(\xi,u) + W(\xi) + \nabla V(\xi) k(\xi, h(x), u) \leq 0$ and $\nabla V(\xi) f(\xi,u) + W(\xi) + \nabla V(\xi) k(\xi, h(x), u) > 0$, using definition (2.8) and noticing that $p(V(\xi)) = 1$ we conclude that (2.12) holds.

Next, we establish the following inequality:

$$(\xi - x)' P\Big(f(\xi, u) + \hat{k}(\xi, h(x), u) - f(x, u)\Big) \leq -c\mu |\xi - x|^2, \text{ for all } (\xi, x, u) \in \Re^n \times \Re^n \times U \text{ with } V(\xi) \leq b \text{ and } V(x) \leq R \tag{2.13}$$

Notice that inequalities (2.1), (2.2) and definitions (2.6), (2.7), (2.8) imply that inequality (2.13) holds for the case $V(\xi) \leq a$. Therefore, we focus on the case $a < V(\xi) \leq b$. Definition (2.7) gives:

$$(\xi - x)' P\Big(f(\xi, u) + \hat{k}(\xi, h(x), u) - f(x, u)\Big) \leq (\xi - x)' P\Big(f(\xi, u) + k(\xi, h(x), u) - f(x, u)\Big) - \frac{\varphi(\xi, h(x), u)}{|\nabla V(\xi)|^2} \nabla V(\xi) P(\xi - x) \tag{2.14}$$

Inequalities (2.2), (2.14) and the fact that $\varphi(\xi, h(x), u) \geq 0$ implies that (2.13) holds if $\nabla V(\xi) P(\xi - x) \geq 0$. Moreover, inequalities (2.2), (2.14) show that (2.13) holds if $\varphi(\xi, h(x), u) = 0$. It remains to consider the case $\nabla V(\xi) P(\xi - x) < 0$ and $\varphi(\xi, h(x), u) > 0$. In this case, definition (2.8) implies $\varphi(\xi, h(x), u) = \nabla V(\xi) f(\xi,u) + W(\xi) + p(V(\xi)) \nabla V(\xi) k(\xi, h(x), u) > 0$. Inequality (2.5) gives:

$$\begin{aligned}
\varphi(\xi, h(x), u)) &= \nabla V(\xi) f(\xi, u) + p(V(\xi)) \nabla V(\xi) k(\xi, h(x), u) + W(\xi) \leq \\
&+ (1 - p(V(\xi))) \nabla V(\xi) f(\xi, u) + (1 - p(V(\xi))) W(\xi) \\
&+ (1 - c) |\nabla V(\xi)|^2 p(V(\xi)) \frac{(\xi - x)' P\big(f(\xi, u) + k(\xi, h(x), u) - f(x, u)\big)}{\nabla V(\xi) P(\xi - x)}
\end{aligned} \tag{2.15}$$

Using (2.15), (2.1) and the fact that $0 \leq p(V(\xi)) \leq 1$, we obtain:

$$\begin{aligned}
-\frac{\varphi(\xi, h(x), u)) \nabla V(\xi) P(\xi - x)}{|\nabla V(\xi)|^2} &\leq \\
-\frac{1 - p(V(\xi))}{|\nabla V(\xi)|^2} \nabla V(\xi) P(\xi - x)\big(\nabla V(\xi) f(\xi, u) + W(\xi)\big) &\\
-(1 - c) p(V(\xi))(\xi - x)' P\big(f(\xi, u) + k(\xi, h(x), u) - f(x, u)\big) &\\
\leq -(1 - c)(\xi - x)' P\big(f(\xi, u) + k(\xi, h(x), u) - f(x, u)\big) &
\end{aligned}$$

Combining (2.2), (2.14) and the above inequality, we conclude that (2.13) holds.

Let arbitrary measurable and locally essentially bounded input $u : \Re_+ \to U$ and arbitrary $(x_0, \xi_0) \in \Re^n \times \Re^n$ and consider the solution $(x(t), \xi(t)) \in \Re^n \times \Re^n$ of (1.1), (1.2) with (2.9), initial condition $(x(0), \xi(0)) = (x_0, \xi_0)$ corresponding to input $u : \Re_+ \to U$. Lemma 2.1 in conjunction with (2.1) and (2.12) implies there exists $T \in C^0(\Re^n; \Re_+)$ such that for every $(x_0, \xi_0) \in \Re^n \times \Re^n$ and for every measurable and locally essentially bounded input $u : \Re_+ \to U$ the solution $(x(t), \xi(t)) \in \Re^n \times \Re^n$ of (1.1), (1.2) with (2.9) with initial condition $(x(0), \xi(0)) = (x_0, \xi_0)$ and corresponding to input $u : \Re_+ \to U$ satisfies:

- $V(x(t)) \leq \max(V(x_0), R)$, $V(\xi(t)) \leq \max(V(\xi_0), b)$ for all $t \geq 0$
- $V(x(t)) \leq R$ for all $t \geq T(x_0)$ and
- $V(\xi(t)) \leq b$ for all $t \geq T(\xi_0)$

Using (2.13) and the absolutely continuous function $Q(t) = (\xi(t) - x(t))' P(\xi(t) - x(t))$, we conclude that

$$|\xi(t) - x(t)| \leq \sqrt{\frac{K_2}{K_1}} \exp(-\sigma(t - t_0)) |\xi(t_0) - x(t_0)|, \text{ for all } t \geq t_0 \qquad (2.16)$$

where $t_0 = \max(T(x_0), T(\xi_0))$, $\sigma := \frac{c\mu}{K_1}$ and $K_2 \geq K_1 > 0$ are constants that satisfy $K_1 |x|^2 \leq x'Px \leq K_2 |x|^2$ for all $x \in \Re^n$. Define:

$$M(x_0, \xi_0) := \sqrt{\frac{K_2}{K_1}} \exp(\sigma \max(T(x_0), T(\xi_0))) \max\{|\xi - x| : V(\xi) \leq \max(V(\xi_0), b), V(x) \leq \max(V(x_0), R)\} \qquad (2.17)$$

Definition (2.17) in conjunction with (2.16) and the fact that $V(x(t)) \leq \max(V(x_0), R)$, $V(\xi(t)) \leq \max(V(\xi_0), b)$ for all $t \geq 0$ implies that (2.10) holds. The proof is complete. ◁

An advantage of the observer design provided by Theorem 2.2 is the fact that the observer can be implemented with sampled measurements. The following result guarantees the design of a global sampled-data exponential observer.

**Theorem 2.4:** *Consider system (1.1), (1.2) under hypotheses (H1-3) and suppose that the following additional hypothesis holds:*

**(H4)** $h(\Re^n) = \Re$ *and there exists a vector* $L \in \Re^n$ *such that* $k(\xi, y, u) = L(h(\xi) - y)$. *Moreover, either* $U \subseteq \Re^m$ *is compact or the mapping* $\nabla h(x) f(x, u)$ *is independent of* $u$.

*Let* $\hat{k} : \Re^n \times \Re \times U \to \Re^n$ *be the locally Lipschitz mapping defined by (2.6), (2.7), (2.8). Then there exists* $M \in C^0(\Re^n \times \Re^n; \Re_+)$ *and* $\sigma, r, \gamma > 0$ *such that for every measurable and locally essentially bounded input* $u : \Re_+ \to U$, *for every locally bounded inputs* $w : \Re_+ \to \Re_+$, $e : \Re_+ \to \Re$ *and for every* $(x_0, \xi_0, w_0) \in \Re^n \times \Re^n \times \Re$ *the solution* $(x(t), \xi(t), w(t)) \in \Re^n \times \Re^n \times \Re$ *of (1.1), (1.2) with*

$$\begin{aligned}
\dot{\xi}(t) &= f(\xi(t), u(t)) + \hat{k}(\xi(t), w(t), u(t)) \\
\dot{w}(t) &= \nabla h(\xi(t)) f(\xi(t), u(t)) \quad , \quad t \in [\tau_i, \tau_{i+1}) \\
w(\tau_{i+1}) &= h(x(\tau_{i+1})) + e(\tau_{i+1}) \\
\tau_{i+1} &= \tau_i + r \exp(-w(\tau_i)) \\
x(0) &= x_0, \xi(0) = \xi_0, w(0) = w_0, \tau_0 = 0
\end{aligned} \qquad (2.18)$$

*satisfies*

$$|\xi(t) - x(t)| \leq M(x_0, \xi_0) \exp(-\sigma t) + \gamma \sup_{0 \leq s \leq t} |e(s)|, \text{ for all } t \geq 0 \qquad (2.19)$$

**Remark 2.5:**
(a) The input $e : \Re_+ \to \Re$ quantifies the effect of measurement errors. It is clear that the sampled-data observer (2.18) satisfies an input-to-output stability property with respect to the measurement error. The input $w : \Re_+ \to \Re_+$ quantifies the effect of perturbations of the sampling schedule. More specifically, the desired inequality (2.19) is guaranteed for every sampling schedule with diameter less or equal than $r > 0$, i.e., for every set of sampling times $\{\tau_i\}_{i=0}^{\infty}$ with $\tau_0 = 0$ and $\sup_{i \geq 0}(\tau_{i+1} - \tau_i) \leq r$ (see also [9,10]). The overall system (1.1), (1.2), (2.18) is a hybrid system with variable sampling partition (see [10]).

**(b)** Hypothesis (H4) holds automatically for nonlinear systems of the form (2.4), where $f_i : \Re^i \to \Re$ ($i=1,...,n$) are smooth mappings.

**(c)** It should be noted that the sampled-data observer (2.18) is similar to the sampled-data observers constructed in [9]. However, the results presented in [9] cannot be used in order to prove Theorem 2.4. The reason is that inequality (3.1) in [9] does not hold for all times (as required in [9]). An analogue of inequality (3.1) in [9] holds after an initial transient period. The transient period is needed so that the state of the original system and the observer state enter an appropriate compact set.

**Proof of Theorem 2.4:** As in the proof of Theorem 2.2, we first notice that for all $(\xi, w, u) \in \Re^n \times \Re \times U$ with $V(\xi) \geq b$ the following inequality holds:

$$\nabla V(\xi)(f(\xi, u) + \hat{k}(\xi, w, u)) \leq -W(\xi) \tag{2.20}$$

Moreover, exactly as in the proof of Theorem 2.2, we show that (2.13) holds. Definitions (2.6), (2.7), (2.8) and hypothesis (H4) imply the existence of a constant $G \geq 0$ such that

$$\left|\hat{k}(\xi, w, u) - \hat{k}(\xi, y, u)\right| \leq G|w - y|, \text{ for all } (\xi, y, w, u) \in \Re^n \times \Re \times \Re \times U \text{ with } V(\xi) \leq b \tag{2.21}$$

$$\left|\nabla h(\xi)f(\xi, u) - \nabla h(x)f(x, u)\right| \leq G|\xi - x|, \text{ for all } (\xi, x, u) \in \Re^n \times \Re^n \times U \text{ with } V(\xi) \leq b \text{ and } V(x) \leq R \tag{2.22}$$

Lemma 2.1 in conjunction with (2.1) and (2.20) implies there exists $T \in C^0(\Re^n; \Re_+)$ such that for every $(x_0, \xi_0, w_0) \in \Re^n \times \Re^n \times \Re$, for every measurable and locally essentially bounded input $u : \Re_+ \to U$ and for every locally bounded inputs $w : \Re_+ \to \Re_+$, $e : \Re_+ \to \Re$, the solution $(x(t), \xi(t), w(t)) \in \Re^n \times \Re^n \times \Re$ of (1.1), (1.2) with (2.18), initial condition $(x(0), \xi(0), w(0)) = (x_0, \xi_0, w_0)$ and corresponding to inputs $u : \Re_+ \to U$, $w : \Re_+ \to \Re_+$, $e : \Re_+ \to \Re$ satisfies:

- $V(x(t)) \leq \max(V(x_0), R)$, $V(\xi(t)) \leq \max(V(\xi_0), b)$ for all $t \geq 0$
- $V(x(t)) \leq R$ for all $t \geq T(x_0)$ and
- $V(\xi(t)) \leq b$ for all $t \geq T(\xi_0)$

It follows from the above estimates that the following inequality holds for all $t \geq 0$:

$$|\xi(t) - x(t)| \leq \max\{|\xi - x| : V(\xi) \leq \max(V(\xi_0), b), V(x) \leq \max(V(x_0), R)\} \tag{2.23}$$

Using (2.13), (2.21) and the absolutely continuous function $Q(t) = (\xi(t) - x(t))' P(\xi(t) - x(t))$, we conclude that

$$|\xi(t) - x(t)| \leq \sqrt{\frac{K_2}{K_1}} \exp(-\sigma(t - t_0))|\xi(t_0) - x(t_0)| + \frac{\sqrt{2}G|P|}{c\mu} \sup_{t_0 \leq s \leq t} (\exp(-\sigma(t - s))|w(s) - y(s)|), \text{ for all } t \geq t_0 \tag{2.24}$$

where $t_0 = \min\{\tau_i : \tau_i \geq \max(T(x_0), T(\xi_0))\}$, $\sigma := \frac{c\mu}{4K_1}$ and $K_2 \geq K_1 > 0$ are constants that satisfy $K_1|x|^2 \leq x'Px \leq K_2|x|^2$ for all $x \in \Re^n$. Finally, notice that for every $t \geq t_0$ the following estimate holds:

$$|w(t) - y(t)| \leq \sup_{0 \leq s \leq t}|e(s)| + rG \sup_{\tau \leq s \leq t}|\xi(s) - x(s)| \tag{2.25}$$

where $\tau = \max\{\tau_i : \tau_i \leq t\}$. Notice that from the inequality $t \leq \tau + r$ and (2.25) we obtain:

$$\exp(\sigma t)|w(t) - y(t)| \leq \exp(\sigma t) \sup_{0 \leq s \leq t}|e(s)| + rG \exp(\sigma(t - \tau)) \sup_{\tau \leq s \leq t}(\exp(\sigma s)|\xi(s) - x(s)|)$$
$$\leq \exp(\sigma t) \sup_{0 \leq s \leq t}|e(s)| + rG \exp(\sigma r) \sup_{\tau \leq s \leq t}(\exp(\sigma s)|\xi(s) - x(s)|) \tag{2.26}$$

Inequalities (2.24) and (2.26) imply the following inequalities for all $t \geq t_0$:

$$\sup_{t_0 \leq s \leq t} \left(\exp(\sigma s)|\xi(s) - x(s)|\right) \leq \sqrt{\frac{K_2}{K_1}} \exp(\sigma t_0)|\xi(t_0) - x(t_0)| + \frac{\sqrt{2}G|P|}{c\mu} \sup_{t_0 \leq s \leq t} \left(\exp(\sigma s)|w(s) - y(s)|\right) \quad (2.27)$$

$$\sup_{t_0 \leq s \leq t} \left(\exp(\sigma s)|w(s) - y(s)|\right) \leq \exp(\sigma t) \sup_{0 \leq s \leq t} |e(s)| + rG \exp(\sigma r) \sup_{t_0 \leq s \leq t} \left(\exp(\sigma s)|\xi(s) - x(s)|\right) \quad (2.28)$$

Finally, we assume that $r > 0$ is selected so that

$$\frac{\sqrt{2}G^2|P|}{c\mu} r \exp(\sigma r) < 1 \quad (2.29)$$

Inequalities (2.27), (2.28) and (2.29) give the following estimate for all $t \geq t_0$:

$$\sup_{t_0 \leq s \leq t} \left(\exp(\sigma s)|\xi(s) - x(s)|\right) \leq \frac{c\mu\gamma}{\sqrt{2}G|P|} \sqrt{\frac{K_2}{K_1}} \exp(\sigma t_0)|\xi(t_0) - x(t_0)| + \gamma \exp(\sigma t) \sup_{0 \leq s \leq t} \left(|e(s)|\right) \quad (2.30)$$

where $\gamma = \dfrac{\sqrt{2}G|P|}{c\mu - \sqrt{2}G^2|P|r\exp(\sigma r)}$. Since $t_0 = \min\{\tau_i : \tau_i \geq \max(T(x_0), T(\xi_0))\}$, it follows that $t_0 \leq \max(T(x_0), T(\xi_0)) + r$. Combining (2.23) and (2.30) we conclude that (2.19) holds with

$$M(x_0, \xi_0) := \frac{c\mu\gamma}{\sqrt{2}G|P|} \sqrt{\frac{K_2}{K_1}} \exp(\sigma r + \sigma \max(T(x_0), T(\xi_0))) \max\left\{|\xi - x| : V(\xi) \leq \max(V(\xi_0), b), V(x) \leq \max(V(x_0), R)\right\}$$

The proof is complete. ◁

## 3. Global Exponential Observers for Systems on Open Sets

Consider the forward complete system:

$$\dot{X} = F(X, u), \quad X \in A, u \in U \quad (3.1)$$

where $A \subseteq \Re^n$ is an open set, $U \subseteq \Re^m$ is a non-empty set, $F : A \times U \to \Re^n$ is a smooth vector field and the output is given by

$$y = H(X) \quad (3.2)$$

where $H : A \to \Re^k$ is a smooth mapping. We assume the knowledge of a smooth injective mapping $\Phi : \Re^n \to A$ with $\Phi(\Re^n) = A$ and $\det(D\Phi(x)) \neq 0$ for all $x \in \Re^n$, where $D\Phi(x) \in \Re^{n \times n}$ is the Jacobian of the mapping $\Phi : \Re^n \to A$, such that system (3.1), (3.2) under the change of coordinates $X = \Phi(x)$ is expressed by (1.1), (1.2), where $f : \Re^n \times U \to \Re^n$ and $h : \Re^n \to \Re^k$ are smooth mappings satisfying

$$D\Phi(x) f(x, u) = F(\Phi(x), u), \text{ for all } (x, u) \in \Re^n \times U \quad (3.3)$$

$$h(x) := H(\Phi(x)), \text{ for all } x \in \Re^n \quad (3.4)$$

where $D\Phi(x) \in \Re^{n \times n}$ is the jacobian of the mapping $\Phi : \Re^n \to A$.

The following hypothesis implies the existence of a "candidate global exponential observer" for system (3.1), (3.2).

**(P1)** *There exist a symmetric and positive definite matrix $P \in \Re^{n \times n}$, a constant $\mu > 0$ and a smooth mapping $k : \Re^n \times H(A) \times U \to \Re^n$ with $k(Z, y, u) = F(Z, u)$ for all $(Z, y, u) \in \Re^n \times H(A) \times U$ with $H(Z) = y$ such that the following inequality holds*

$$(Z - X)'P\big(k(Z, H(X), u) - F(X, u)\big) \leq -\mu |Z - X|^2, \text{ for all } u \in U, (Z, X) \in \Re^n \times A \quad (3.5)$$

Indeed, hypothesis (P1) guarantees that for every $(X(0), Z(0)) \in A \times \Re^n$ and for every measurable and locally essentially bounded input $u : \Re_+ \to U$ the solution of system (3.1), (3.2) with

$$\dot{Z} = k(Z, y, u) \quad (3.6)$$

will satisfy an estimate of the form $|Z(t) - X(t)| \leq M \exp(-\sigma t)|Z(0) - X(0)|$, for all $t \geq 0$ for appropriate constants $M, \sigma > 0$. However, system (3.6) is not necessarily an observer, since we cannot guarantee that $Z(t) \in A$ for all $t \geq 0$.

In order to state the problem in a different way, it is convenient to use the change of coordinates $Z = \Phi(z)$ for the "observer" (3.6):

$$\dot{z} = \widetilde{k}(z, y, u), \quad z \in \Re^n \quad (3.7)$$

where $\widetilde{k}(z, y, u) := (D\Phi(z))^{-1} k(\Phi(z), y, u)$ is a smooth mapping. Now, the problem can be stated as follows:

"Although system (1.1), (1.2) is forward complete,
system (1.1), (1.2) with (3.7) is not necessarily forward complete"

Since system (1.1), (1.2) is forward complete, the results in [3] guarantee the existence of a radially unbounded (but not necessarily positive definite) function $W \in C^2(\Re^n; [1, +\infty))$, a continuous function $K : U \to \Re_+$ and a constant $R \geq 0$ such that

$$\nabla W(x) f(x, u) \leq K(u) W(x), \text{ for all } (x, u) \in \Re^n \times U \text{ with } W(x) \geq R \quad (3.8)$$

The problem that we consider in this section is the problem of existence/design of an observer with state $Z \in A$ which guarantees global exponential convergence in the original $(X, Z)$ coordinates based on the knowledge of the function $W$ and the "candidate observer" (3.6). Our main result guarantees that under some additional assumptions the existence/design problem of the global exponential observer is solvable.

**Theorem 3.1:** *Suppose that there exist constants $a > R$, $\varepsilon \in (0, 1)$, a symmetric and positive semidefinite matrix $Q(z) = \{q_{i,j}(z) : i, j = 1, ..., n\} \in \Re^{n \times n}$ with $q_{i,j} \in C^1(\Re^n; \Re)$ and a continuous function $c : H(A) \times U \to [1, +\infty)$ such that the following inequality holds for all $(z, x, u) \in \Re^n \times \Re^n \times U$ with $W(z) \geq a$ and $\nabla W(z) Q(z) (D\Phi(z))' P(\Phi(z) - \Phi(x)) < 0$*

$$\nabla W(z) \widetilde{k}(z, y, u) \leq c(y, u) W(z) + (1 - \varepsilon) \nabla W(z) Q(z) (\nabla W(z))' \frac{(\Phi(z) - \Phi(x))' P \big(D\Phi(z) \widetilde{k}(z, y, u) - D\Phi(x) f(x, u)\big)}{\nabla W(z) Q(z) (D\Phi(z))' P(\Phi(z) - \Phi(x))} \quad (3.9)$$

*Define $C := \big\{ z \in \Re^n : \nabla W(z) Q(z) (\nabla W(z))' = 0, W(z) \geq a \big\}$. Suppose that either the set $C$ is empty or that $\nabla W(z) \widetilde{k}(z, y, u) < c(y, u) W(z)$ for all $(z, y, u) \in C \times H(A) \times U$.*

*Define the locally Lipschitz mapping* $G: A \times H(A) \times U \to \Re^n$ *by*

$$G(Z, y, u) := k(Z, y, u) - \lambda(\Phi^{-1}(Z), y, u) D\Phi(\Phi^{-1}(Z)) Q(\Phi^{-1}(Z)) (\nabla W(\Phi^{-1}(Z)))',$$
$$\text{for all } (Z, y, u) \in A \times H(A) \times U \quad (3.10)$$

*where* $\lambda : \Re^n \times H(A) \times U \to \Re_+$ *is defined by*

$$\lambda(z, y, u) := \begin{cases} \dfrac{p(W(z)) \nabla W(z) \widetilde{k}(z, y, u) - c(y, u) W(z)}{\nabla W(z) Q(z) (\nabla W(z))'} & \text{if } p(W(z)) \nabla W(z) \widetilde{k}(z, y, u) > c(y, u) W(z) \\ 0 & \text{if } p(W(z)) \nabla W(z) \widetilde{k}(z, y, u) \leq c(y, u) W(z) \end{cases},$$
$$\text{for all } (z, y, u) \in \Re^n \times H(A) \times U \quad (3.11)$$

*and* $p: \Re_+ \to [0,1]$ *is an arbitrary locally Lipschitz function that satisfies* $p(s) = 1$ *for all* $s \geq a+1$ *and* $p(s) = 0$ *for all* $s \leq a$. *Then there exists* $M > 0$ *such that for every measurable and locally essentially bounded input* $u: \Re_+ \to U$ *and for every* $(X_0, Z_0) \in A \times A$ *the solution* $(X(t), Z(t)) \in A \times A$ *of (3.1), (3.2) with*

$$\dot{Z} = G(Z, y, u), Z \in A \quad (3.12)$$

*initial condition* $(X(0), Z(0)) = (X_0, Z_0)$ *corresponding to input* $u: \Re_+ \to U$ *satisfies*

$$|Z(t) - X(t)| \leq M |X_0 - Z_0| \exp\left(-\frac{\varepsilon \mu}{2} t\right), \text{ for all } t \geq 0 \quad (3.13)$$

**Proof:** Inequality (3.5), in conjunction with definitions (3.3), (3.4) and definition $\widetilde{k}(z, y, u) := (D\Phi(z))^{-1} k(\Phi(z), y, u)$, implies that the following inequality holds:

$$(\Phi(z) - \Phi(x))' P\left(D\Phi(z) \widetilde{k}(z, h(x), u) - D\Phi(x) f(x, u)\right) \leq -\mu |\Phi(z) - \Phi(x)|^2, \text{ for all } u \in U, (z, x) \in \Re^n \times \Re^n \quad (3.14)$$

We next evaluate the quantity $\dot{W}(z, x, u) := \nabla W(z) \hat{k}(z, h(x), u)$ for all $u \in U$, $(z, x) \in \Re^n \times \Re^n$ with $W(z) \geq a+1$, where $\hat{k}(z, y, u) := \widetilde{k}(z, y, u) - \lambda(z, y, u) Q(z) (\nabla W(z))'$. We get:

$$\dot{W}(z, x, u) = \nabla W(z) \widetilde{k}(z, h(x), u) - \lambda(z, h(x), u) \nabla W(z) Q(z) (\nabla W(z))'$$

Since $W(z) \geq a+1$ and $p: \Re_+ \to [0,1]$ is an arbitrary locally Lipschitz function that satisfies $p(s) = 1$ for all $s \geq a+1$, we obtain $p(W(z)) = 1$. By distinguishing the cases $p(W(z)) \nabla W(z) \widetilde{k}(z, y, u) > c(y, u) W(z)$ and $p(W(z)) \nabla W(z) \widetilde{k}(z, y, u) \leq c(y, u) W(z)$ we may conclude that the following inequality holds for all $u \in U$, $(z, x) \in \Re^n \times \Re^n$ with $W(z) \geq a+1$:

$$\dot{W}(z, x, u) \leq c(h(x), u) W(z) \quad (3.15)$$

We next claim that the following system is forward complete:

$$\begin{aligned} \dot{x} &= f(x, u) \\ \dot{z} &= \widetilde{k}(z, h(x), u) - \lambda(z, h(x), u) Q(z) (\nabla W(z))' \\ z &\in \Re^n, x \in \Re^n, u \in U \end{aligned} \quad (3.16)$$

The claim is proved at the Appendix. Since system (3.16) is forward complete and using the change of coordinates $X = \Phi(x)$, $Z = \Phi(z)$, we conclude that the solution $(X(t), Z(t))$ of (3.1), (3.2) and (3.6) starting from arbitrary initial condition $(X(0), Z(0)) \in A \times A$ and corresponding to arbitrary measurable and locally essentially bounded input $u : \Re_+ \to U$ exists for all $t \geq 0$ and satisfies $(X(t), Z(t)) \in A \times A$ for all $t \geq 0$.

We next evaluate the quantity $(\Phi(z) - \Phi(x))' P \big( D\Phi(z) \hat{k}(z, h(x), u) - D\Phi(x) f(x, u) \big)$, where $\hat{k}(z, y, u) := \tilde{k}(z, y, u) - \lambda(z, y, u) Q(z) (\nabla W(z))'$. Consider the cases:

1) $W(z) \leq a$. In this case, $\lambda(z, y, u) = 0$ and $\hat{k}(z, y, u) = \tilde{k}(z, y, u)$. Therefore, inequality (3.14) implies that $(\Phi(z) - \Phi(x))' P \big( D\Phi(z) \hat{k}(z, h(x), u) - D\Phi(x) f(x, u) \big) \leq -\mu |\Phi(z) - \Phi(x)|^2$.

2) $W(z) > a$. In this case $\hat{k}(z, y, u) = \tilde{k}(z, y, u) - \lambda(z, y, u) Q(z) (\nabla W(z))'$ and we get:

$$\begin{aligned}
&(\Phi(z) - \Phi(x))' P \big( D\Phi(z) \hat{k}(z, h(x), u) - D\Phi(x) f(x, u) \big) = \\
&= (\Phi(z) - \Phi(x))' P \Big( D\Phi(z) \tilde{k}(z, h(x), u) - \lambda(z, h(x), u) D\Phi(z) Q(z) (\nabla W(z))' - D\Phi(x) f(x, u) \Big) = \\
&= (\Phi(z) - \Phi(x))' P \big( D\Phi(z) \tilde{k}(z, h(x), u) - D\Phi(x) f(x, u) \big) - \lambda(z, h(x), u) \nabla W(z) Q(z) (D\Phi(z))' P(\Phi(z) - \Phi(x))
\end{aligned} \quad (3.17)$$

Notice that $\lambda(z, y, u) \geq 0$. If $\nabla W(z) Q(z) (D\Phi(z))' P(\Phi(z) - \Phi(x)) \geq 0$ then (3.17) implies $(\Phi(z) - \Phi(x))' P \big( D\Phi(z) \hat{k}(z, h(x), u) - D\Phi(x) f(x, u) \big) \leq -\mu |\Phi(z) - \Phi(x)|^2$.

If $\nabla W(z) Q(z) (D\Phi(z))' P(\Phi(z) - \Phi(x)) < 0$ and $\lambda(z, y, u) = 0$ then (3.17) implies $(\Phi(z) - \Phi(x))' P \big( D\Phi(z) \hat{k}(z, h(x), u) - D\Phi(x) f(x, u) \big) \leq -\mu |\Phi(z) - \Phi(x)|^2$.

Thus the only case that remains to be considered is the case $\nabla W(z) Q(z) (D\Phi(z))' P(\Phi(z) - \Phi(x)) < 0$ and $\lambda(z, y, u) > 0$. In this case, we have $\lambda(z, y, u) = \dfrac{p(W(z)) \nabla W(z) \tilde{k}(z, y, u) - c(y, u) W(z)}{\nabla W(z) Q(z) (\nabla W(z))'}$ and $p(W(z)) \nabla W(z) \tilde{k}(z, y, u) > c(y, u) W(z)$. Inequality (3.9) gives:

$$\lambda(z, y, u) \leq \frac{(p(W(z)) - 1) c(y, u) W(z)}{\nabla W(z) Q(z) (\nabla W(z))'} + (1 - \varepsilon) p(W(z)) \frac{(\Phi(z) - \Phi(x))' P \big( D\Phi(z) \tilde{k}(z, y, u) - D\Phi(x) f(x, u) \big)}{\nabla W(z) Q(z) (D\Phi(z))' P(\Phi(z) - \Phi(x))}$$

The above inequality in conjunction with (3.14) and (3.17) implies:

$$(\Phi(z) - \Phi(x))' P \big( D\Phi(z) \hat{k}(z, h(x), u) - D\Phi(x) f(x, u) \big) \leq -\mu \varepsilon |\Phi(z) - \Phi(x)|^2 \quad (3.18)$$

Since $\varepsilon \in (0, 1)$, we conclude that the above inequality holds for all $u \in U$, $(z, x) \in \Re^n \times \Re^n$. using the change of coordinates $X = \Phi(x)$, $Z = \Phi(z)$ and the differential inequality (3.18), we conclude that the solution $(X(t), Z(t))$ of (3.1), (3.2) and (3.6) starting from arbitrary initial condition $(X(0), Z(0)) \in A \times A$ and corresponding to arbitrary measurable and locally essentially bounded input $u : \Re_+ \to U$ satisfies for almost all $t \geq 0$:

$$\dot{V}(t) \leq -\mu \varepsilon V(t)$$

where $V(t) = (Z(t) - X(t))' P (Z(t) - X(t))$. The existence of a constant $M > 0$ satisfying (3.13) is a direct consequence of the above differential inequality. The proof is complete. ◁

## 4. Examples

In the present section, we will apply the results of the previous Sections to two specific examples. The first example is a system with monotone nonlinearities, and we will apply the results of Theorems 2.2 and 2.4. The second example is a chemostat model, with positive state variables, and we apply the construction of Section 3.

**Example 4.1:** Consider the nonlinear planar system

$$\dot{x}_1 = -x_1^3 + x_2$$
$$\dot{x}_2 = -x_2^3 + u$$
$$y = x_1 \qquad (4.1)$$
$$(x_1, x_2)' \in \Re^2, u \in U = [-1,1]$$

System (4.1) is a system of the form (2.4) with monotone nonlinearities. The nonlinearities are not globally Lipschitz; however a continuous global exponential observer can be designed using the methodology proposed in [1,7]. Here, we will design a continuous global exponential observer using Theorem 2.2 and we will show that we can also design a robust global exponential sampled-data observer using Theorem 2.4.

We will show next that hypotheses (H1-4) hold for system (4.1). First notice that hypothesis (H1) holds with $V(x) = \frac{1}{2}x_1^2 + \frac{1}{2}x_2^2$. Indeed, notice that the inequalities $x_1 x_2 \leq \frac{1}{2}x_1^2 + \frac{1}{2}x_2^2$, $x_2 u \leq \frac{1}{2}x_2^2 + \frac{1}{2}u^2$, $u^2 \leq 1$, $\frac{1}{2}x_1^2 \leq \frac{1}{4}x_1^4 + \frac{1}{4}$, $x_2^2 \leq \frac{1}{2}x_2^4 + \frac{1}{2}$ and $V^2(x) \leq \frac{1}{2}x_1^4 + \frac{1}{2}x_2^4$ give us:

$$\nabla V(x) f(x,u) = x_1 x_2 + x_2 u - x_1^4 - x_2^4 \leq \frac{1}{2}x_1^2 + x_2^2 + \frac{1}{2}u^2 - x_1^4 - x_2^4$$
$$\leq \frac{3}{4} + \frac{1}{2}u^2 - \frac{3}{4}x_1^4 - \frac{1}{2}x_2^4 \leq \frac{5}{4} - \frac{1}{2}x_1^4 - \frac{1}{2}x_2^4 \leq \frac{5}{4} - V^2(x)$$

where $f(x,u) := \left(-x_1^3 + x_2, \ -x_2^3 + u\right)'$. The above inequality shows that inequality (2.1) holds with $W(x) = \frac{1}{2}V^2(x)$ and $R = \frac{\sqrt{10}}{2}$.

Next, we show that hypothesis (H2) holds. Let $P = \frac{1}{2}\begin{bmatrix} 1 & -p \\ -p & q \end{bmatrix}$ with $q > p^2$, $p > 0$ and $k(\xi, y, u) := \begin{bmatrix} L_1 \\ L_2 \end{bmatrix}(\xi_1 - y)$, where $p, q, L_1, L_2$ are constants to be selected. We get:

$$2(\xi - x)' P(f(\xi, u) + k(\xi, y, u) - f(x, u)) =$$
$$= (1 - pL_1 + qL_2)(\xi_1 - x_1)(\xi_2 - x_2) + (L_1 - pL_2)(\xi_1 - x_1)^2 - p(\xi_2 - x_2)^2$$
$$+ p(\xi_1 - x_1)(\xi_2^3 - x_2^3) - q(\xi_2 - x_2)(\xi_2^3 - x_2^3) - (\xi_1 - x_1)(\xi_1^3 - x_1^3) + p(\xi_1^3 - x_1^3)(\xi_2 - x_2)$$

Using the inequalities $-(\xi_1 - x_1)(\xi_1^3 - x_1^3) \leq 0$, $p(\xi_1 - x_1)(\xi_2^3 - x_2^3) \leq q(\xi_2 - x_2)(\xi_2^3 - x_2^3) + \frac{p}{4q}(\xi_2^2 + \xi_2 x_2 + x_2^2)(\xi_1 - x_1)^2$, $p(\xi_1^3 - x_1^3)(\xi_2 - x_2) \leq \frac{p}{2}(\xi_2 - x_2)^2 + \frac{p}{2}(\xi_1 - x_1)^2(\xi_1^2 + \xi_1 x_1 + x_1^2)^2$, we obtain:

$$2(\xi - x)' P(f(\xi, u) + k(\xi, y, u) - f(x, u)) \leq$$
$$\leq (1 - pL_1 + qL_2)(\xi_1 - x_1)(\xi_2 - x_2) - \frac{p}{2}(\xi_2 - x_2)^2$$
$$+ \left(L_1 - pL_2 + \frac{p}{4q}(\xi_2^2 + \xi_2 x_2 + x_2^2) + \frac{p}{2}(\xi_1^2 + \xi_1 x_1 + x_1^2)^2\right)(\xi_1 - x_1)^2$$

If we further assume that $V(x) \leq R$ and $V(\xi) \leq b$, where $b > R$ is a constant to be selected, we obtain:

$$2(\xi - x)' P(f(\xi,u) + k(\xi,y,u) - f(x,u)) \leq$$
$$\leq (1 - pL_1 + qL_2)(\xi_1 - x_1)(\xi_2 - x_2) - \frac{p}{2}(\xi_2 - x_2)^2 + \left(L_1 - pL_2 + \frac{3pb}{2q} + 18pb^2\right)(\xi_1 - x_1)^2$$

The selection

$$L_1 = -p\frac{q + 2 + 3b + 36qb^2}{2(q - p^2)}, \quad L_2 = -\frac{1}{q}\left(p^2 \frac{q + 2 + 3b + 36qb^2}{2(q - p^2)} + 1\right) \quad (4.2)$$

in conjunction with the above inequality implies that (2.2) holds with $\mu = \frac{p}{4}$. Therefore, we have showed that hypothesis (H2) holds with arbitrary $p, q, b$ satisfying $q > p^2$, $p > 0$ and $b > R$. Moreover, hypothesis (H4) holds as well (notice that $h(\xi) = \xi_1$).

Finally, we show that hypothesis (H3) holds. More specifically, we will show that there exist constants $R \leq a < b$ such that (2.11) holds for all $u \in U = [-1,1]$, $\xi, x \in \Re^2$ with $a < V(\xi) \leq b$ and $V(x) \leq R$ with $c = \frac{1}{2}$ and $W(x) = \frac{1}{2}V^2(x)$. Inequality (2.11) is equivalent to the following inequality:

$$L_1(\xi_1 - x_1)\xi_1 + L_2(\xi_1 - x_1)\xi_2 + \xi_1\xi_2 - \xi_1^4 + \xi_2 u - \xi_2^4 \leq -\frac{1}{2}V^2(\xi) + \frac{p}{8|P|}|\xi||\xi - x| \quad (4.3)$$

Using the inequalities $\xi_1\xi_2 \leq \frac{1}{2}\xi_1^2 + \frac{1}{2}\xi_2^2$, $\xi_2 u \leq \frac{1}{2}\xi_2^2 + \frac{1}{2}u^2$, $u^2 \leq 1$, $\frac{1}{2}\xi_1^2 \leq \frac{1}{4}\xi_1^4 + \frac{1}{4}$, $\xi_2^2 \leq \frac{1}{2}\xi_2^4 + \frac{1}{2}$ and $V^2(\xi) \leq \frac{1}{2}\xi_1^4 + \frac{1}{2}\xi_2^4$, $|L_1\xi_1| + |L_2\xi_2| \leq \sqrt{2}\max(|L_1|,|L_2|)|\xi|$, we conclude that (4.3) holds provided that the following (more demanding) inequality holds for all $\xi, x \in \Re^2$ with $a < V(\xi) \leq b$ and $V(x) \leq R$:

$$|\xi_1 - x_1|\sqrt{2}\max(|L_1|,|L_2|)|\xi| + \frac{5}{4} \leq \frac{1}{2}V^2(\xi) \quad (4.4)$$

Since (4.4) must hold for all $\xi, x \in \Re^2$ with $a < V(\xi) \leq b$ and $V(x) \leq R < b$, we conclude that (4.4) holds automatically provided the following inequality holds:

$$\max(|L_1|,|L_2|) \leq \frac{2a^2 - 5}{16\sqrt{2}b} \quad (4.5)$$

Definitions (4.2) imply that inequality (4.5) holds provided that $b > a > R = \frac{\sqrt{10}}{2}$ and $p, q > 0$ are selected so that:

$$q > \frac{16\sqrt{2}b}{2a^2 - 5} \text{ and } p \leq \min\left(q\frac{2a^2 - 5}{16\sqrt{2}b(q + 2 + 3b + 36qb^2)}, \sqrt{\frac{q}{2}}, q - \frac{16\sqrt{2}b}{2a^2 - 5}\right) \quad (4.6)$$

For example, all inequalities hold for $a = 2, b = 3$, $q = 32$, $p = \frac{\sqrt{2}}{10411}$. Therefore, we are in a position to define the mapping:

$$\hat{k}(\xi, y, u) := \begin{bmatrix} L_1 \\ L_2 \end{bmatrix} (\xi_1 - y), \text{ for all } (\xi, y, u) \in \Re^2 \times \Re \times [-1,1] \text{ with } |\xi|^2 \leq \sqrt{10} \quad (4.7)$$

$$\hat{k}(\xi, y, u) := \begin{bmatrix} L_1(\xi_1 - y) - \dfrac{\varphi(\xi, y, u)}{|\xi|^2} \xi_1 \\ L_2(\xi_1 - y) - \dfrac{\varphi(\xi, y, u)}{|\xi|^2} \xi_2 \end{bmatrix}, \text{ for all } (\xi, y, u) \in \Re^2 \times \Re \times [-1,1] \text{ with } |\xi|^2 > \sqrt{10} \quad (4.8)$$

where $\varphi : \Re^2 \times \Re \times [-1,1] \to \Re_+$ is defined by

$$\varphi(\xi, y, u) := \max\left( 0, -\xi_1^4 + (\xi_1 + u)\xi_2 - \xi_2^4 + \frac{1}{2}V^2(\xi) + (\xi_1 - y)g(V(\xi))(L_1\xi_1 + L_2\xi_2) \right) \quad (4.9)$$

and $g : \Re_+ \to [0,1]$ is an arbitrary locally Lipschitz function that satisfies $g(s) = 1$ for all $s \geq 3$ and $g(s) = 0$ for all $s \leq 2$. Theorem 2.4 guarantees the existence of $r > 0$ such that the system

$$\begin{aligned}
\dot{\xi}_1(t) &= -\xi_1^3(t) + \xi_2(t) + \hat{k}_1(\xi(t), w(t), u(t)) \\
\dot{\xi}(t) &= -\xi_2^3(t) + u(t) + \hat{k}_2(\xi(t), w(t), u(t)) \\
\dot{w}(t) &= -\xi_1^3(t) + \xi_2(t) \quad , \quad t \in [\tau_i, \tau_{i+1}) \\
w(\tau_{i+1}) &= x_1(\tau_{i+1}) + e(\tau_{i+1}) \\
\tau_{i+1} &= \tau_i + r \exp(-w(\tau_i))
\end{aligned} \quad (4.10)$$

is a robust global sampled-data exponential observer. ◁

**Example 4.2:** Consider the chemostat model ([15]):

$$\begin{aligned}
\dot{X} &= X(\mu(S) - D - b) \\
\dot{S} &= D(S_i - S) - K\mu(S)X \\
X &\in (0, +\infty), S \in (0, +\infty)
\end{aligned} \quad (4.11)$$

with output $y = \mu(S)X$ and inputs $S_i : \Re_+ \to [\theta, +\infty)$, $D : \Re_+ \to [\theta, +\infty)$, where $K, \theta > 0$, $b \geq 0$ are constants and $\mu : \Re_+ \to [0, \mu_{\max}]$ is a locally Lipschitz bounded function with $\mu(0) = 0$ and $\mu(S) > 0$ for all $S > 0$. Physically, the system states $X$ and $S$ represent the biomass concentration and substrate concentration respectively, both positive quantities, and $A = \text{int}(\Re_+^2)$ is a positively invariant open set that contains the physically meaningful trajectories of the system. The term $\mu(S)X$ represents the growth rate of microorganisms and it is a measurable quantity in bioreactors with a gaseous product, like anaerobic digesters, where the biogas production rate is proportional to the microbial growth rate ([6,12]).

The following "candidate observer":

$$\begin{aligned}
\dot{Z}_1 &= -(D+b)Z_1 + y \\
\dot{Z}_2 &= D(S_i - Z_2) - Ky \\
Z(t) &= (Z_1(t), Z_2(t))
\end{aligned} \quad (4.12)$$

satisfies hypothesis (P1) with $P = \dfrac{1}{2}I$, $\mu = \theta > 0$. However, it is clear that the "candidate observer" (4.12) does not satisfy the requirement $Z(t) \in A = \text{int}(\Re_+^2)$ for all initial conditions and all times. Here, we will apply the results of Section 3 using the smooth injective mapping $\Phi : \Re^n \to A$ defined by

$$\Phi(x) = \begin{bmatrix} e^{x_1} \\ e^{x_2} \end{bmatrix} \quad (4.13)$$

We will construct a global exponential observer for system (4.11) under the assumption that there exists $S^* > 0$ such that

$$(S - S')(\mu(S) - \mu(S')) \geq 0, \text{ for all } S, S' \in [0, S^*] \quad (4.14)$$

i.e., we will assume that $\mu$ is non-decreasing on the interval $[0, S^*]$.

Indeed, the smooth injective mapping $\Phi : \Re^n \to A$ defined by (4.14) allows us to determine the vector fields
$f(x,u) = \begin{bmatrix} \mu(e^{x_2}) - u_1 - b \\ u_1(u_2 e^{-x_2} - 1) - K\mu(e^{x_2})e^{x_1 - x_2} \end{bmatrix}$, $k(Z, y, u) = \begin{bmatrix} -(u_1 + b)Z_1 + y \\ u_1(u_2 - Z_2) - Ky \end{bmatrix}$, $\widetilde{k}(z, y, u) = \begin{bmatrix} -(u_1 + b) + ye^{-z_1} \\ u_1(u_2 e^{-z_2} - 1) - Kye^{-z_2} \end{bmatrix}$ so
that all equations (3.1)-(3.7) hold with $h(x) := \mu(e^{x_2})e^{x_1}$, $u = (u_1, u_2)$, $u_1 = D$, $u_2 = S_i$ and $U = [\theta, +\infty) \times [\theta, +\infty)$. A radially unbounded (but not positive definite) function $W \in C^2(\Re^n; [4, +\infty))$ may be defined by:

$$W(x) = e^{x_1} + 3e^{-x_1} + e^{2x_1} + e^{x_2} + e^{-x_2} \quad (4.15)$$

Since $\mu(e^{x_2}) \in (0, \mu_{\max}]$, $u_1 \geq \theta$, $u_2 \geq \theta$, we get for all $(x, u) \in \Re^2 \times U$:

$$\nabla W(x)f(x,u) \leq (\mu_{\max} - \theta - b)e^{x_1} + 2(\mu_{\max} - \theta - b)e^{2x_1} + 3(u_1 + b)e^{-x_1} + u_1 u_2$$
$$+ u_1 e^{-x_2} + K\mu(e^{x_2})e^{x_1 - 2x_2} - \theta^2 e^{-2x_2} \quad (4.16)$$

Since $\mu : \Re_+ \to [0, \mu_{\max}]$ is a locally Lipschitz bounded function with $\mu(0) = 0$ and $\mu(S) > 0$ for all $S > 0$, there exists a constant $\gamma > 0$ such that $\mu(S) \leq \gamma S$ for all $S > 0$. Therefore (4.16) implies that the following differential inequality holds for all $(x, u) \in \Re^2 \times U$:

$$\nabla W(x)f(x,u) \leq (\mu_{\max} - \theta - b)e^{x_1} + 2(\mu_{\max} - \theta - b)e^{2x_1} + 3(u_1 + b)e^{-x_1} + u_1 u_2$$
$$+ u_1 e^{-x_2} + K\gamma e^{x_1 - x_2} - \theta^2 e^{-2x_2}$$

Using the inequalities $K\gamma e^{x_1 - x_2} \leq \theta^2 e^{-2x_2} + \left(\dfrac{K\gamma}{2\theta}\right)^2 e^{2x_1}$, $e^{x_2} + e^{-x_2} \geq 2$ in conjunction with the above differential inequality, we conclude that (3.8) holds with $K(u) := \max\left\{1, \mu_{\max} - \theta - b, 2(\mu_{\max} - \theta - b) + \left(\dfrac{K\gamma}{2\theta}\right)^2, u_1 + b, \dfrac{1}{2}u_1 u_2 + u_1\right\}$ and arbitrary constant $R \geq 0$. Finally, we evaluate the quantity $\nabla W(z)\widetilde{k}(z, y, u)$ for all $(z, u) \in \Re^2 \times U$ and $y > 0$:

$$\nabla W(z)\widetilde{k}(z, y, u) \leq y(1 + e^{z_1}) + 3(u_1 + b)e^{-z_1} + u_1 u_2 + u_1 e^{-z_2} - \theta^2 e^{-2z_2} + Kye^{-2z_2}$$

Using the inequalities $e^{z_2} + e^{-z_2} \geq 2$, $e^{z_1} + e^{-z_1} \geq 2$ and the above inequality we get for all $(z, u) \in \Re^2 \times U$ and $y > 0$:

$$\nabla W(z)\widetilde{k}(z, y, u) \leq \max\left\{\dfrac{3}{2}y, u_1 + b, u_1 + \dfrac{1}{2}u_1 u_2\right\} W(z) + (Ky - \theta^2)e^{-2z_2} \quad (4.17)$$

$$\nabla W(z)\widetilde{k}(z, y, u) \leq y(1 + e^{z_1}) + 3(u_1 + b)e^{-z_1} + u_1 u_2 + u_1 e^{-z_2} + K\mu(e^{z_2})e^{z_1 - 2z_2} - \theta^2 e^{-2z_2} + K(y - \mu(e^{z_2})e^{z_1})e^{-2z_2} \quad (4.18)$$

Again using the fact that there exists a constant $\gamma > 0$ such that $\mu(S) \leq \gamma S$ for all $S > 0$ and the inequality $K\gamma e^{z_1 - z_2} \leq \theta^2 e^{-2z_2} + \left(\frac{K\gamma}{2\theta}\right)^2 e^{2z_1}$, we obtain from (4.17) and (4.18) for all $(z, u) \in \Re^2 \times U$ and $y > 0$:

$$\nabla W(z) \widetilde{k}(z, y, u) \leq \max\left\{\frac{3y}{2}, u_1 + b, u_1 + \frac{u_1 u_2}{2}, \left(\frac{K\gamma}{2\theta}\right)^2\right\} W(z) + \min\left(Ky - \theta^2, K\left(y - \mu(e^{z_2})e^{z_1}\right)\right) e^{-2z_2} \quad (4.19)$$

Define $c(y, u) := \max\left\{\frac{3}{2}y, u_1 + b, u_1 + \frac{1}{2}u_1 u_2, \left(\frac{K\gamma}{2\theta}\right)^2\right\} + \frac{Ke^{2p}}{4} y + 1$, $Q(z) = \begin{bmatrix} 0 & 0 \\ 0 & 1 \end{bmatrix}$, where $p > 0$ satisfies $-p \leq \ln(S^*)$ and is yet to be selected. Inequality (4.19) implies that $\nabla W(z) \widetilde{k}(z, y, u) < c(y, u) W(z)$ and (3.9) holds with arbitrary $\varepsilon \in (0,1)$ for all $(z, u) \in \Re^2 \times U$ and $y > 0$ with $z_2 \geq -p$. It follows from (4.19) that (3.9) holds provided that there exists constant $\varepsilon \in (0,1)$ such that the following inequality holds for all $(z, x) \in \Re^2 \times \Re^2$ with $x_2 < z_2 < -p$:

$$K\left(\mu(e^{x_2})e^{x_1} - \mu(e^{z_2})e^{z_1}\right) \leq 2\theta(1 - \varepsilon)\left(1 - e^{-2p}\right)\left[\frac{(e^{z_1} - e^{x_1})^2}{|e^{z_2} - e^{x_2}|} + |e^{z_2} - e^{x_2}|\right] \quad (4.20)$$

Using the fact that $\mu$ is non-decreasing on the interval $[0, S^*]$, we conclude that (4.20) holds for all $(z, x) \in \Re^2 \times \Re^2$ with $x_2 < z_2 < -p$, provided that the following inequality holds:

$$K|e^{x_1} - e^{z_1}|\mu(e^{-p}) \leq 2\theta(1 - \varepsilon)\left(1 - e^{-2p}\right)\left[\frac{(e^{z_1} - e^{x_1})^2}{|e^{z_2} - e^{x_2}|} + |e^{z_2} - e^{x_2}|\right] \quad (4.21)$$

Since $\mu(0) = 0$, continuity of $\mu$ implies that there exists sufficiently large $p > 0$ such that $K\mu(e^{-p}) \leq 4\theta(1 - \varepsilon)\left(1 - e^{-2p}\right)$, which directly implies inequality (4.21). Therefore, we conclude that (3.9) holds with arbitrary $a \geq 0$. The global exponential observer will be given by the equations:

$$\begin{aligned} \dot{Z}_1 &= (u_1 + b)Z_1 + y \\ \dot{Z}_2 &= u_1(u_2 - Z_2) - Ky + \lambda(\Phi^{-1}(Z), y, u)\left(1 - Z_2^2\right) \end{aligned} \quad (4.22)$$

where $\lambda$ is defined by (3.11) and $p : \Re_+ \to [0,1]$ is an arbitrary locally Lipschitz function that satisfies $p(s) = 1$ for all $s \geq a + 1$ and $p(s) = 0$ for all $s \leq a$. ◁

## 5. Concluding Remarks

This work developed sufficient conditions for the existence of global exponential observers for two classes of nonlinear systems. The first is the class of systems with a globally asymptotically stable compact set. The second is the class of systems that evolve on an open proper subset of $\Re^n$. In both cases, the construction starts with a "candidate observer", which is subsequently modified by adding a correction term, in order to satisfy appropriate Lyapunov inequalities. In the first class of systems, the "candidate observer" is a local observer over a certain compact set, whereas the correction term forces the trajectory to enter the compact set in finite time. In the second class of systems, the "candidate observer" does not guarantee that the observer trajectories lie within the open set, but this is accomplished through an appropriate correction term. The design of the correction term is performed after transforming the system through an appropriate smooth injective map that maps the open set onto $\Re^n$.

The derived continuous-time observer can lead to the construction of a robust global sampled-data exponential observer. The ideas developed to handle the second class of systems could find potential use in the context of transformation-based observers, relaxing the requirement of a diffeomorphism of $\Re^n$ onto $\Re^n$, allowing the image of the inverse map to be an open subset of $\Re^n$.

# Appendix

**Proof of Lemma 2.1:** The fact that the compact set $S = \{x \in \Re^n : V(x) \leq R\}$ is positively invariant for every measurable and locally essentially bounded input $u : \Re_+ \to U$ is a direct consequence of differential inequality (2.1). Next, we consider the solution $x(t) \in \Re^n$ of (1.1) with arbitrary initial condition $x(0) = x_0 \notin S$ (i.e., $V(x_0) > R$) and corresponding to arbitrary input $u : \Re_+ \to U$.

There exists $t_{\max} \in (0,+\infty]$ such that the solution is defined on $[0, t_{\max})$ (and cannot be extended if $t_{\max} < +\infty$). Define $A := \{t \in [0, t_{\max}) : V(x(t)) \leq R\}$. We will show next that $A \neq \emptyset$. Suppose that $A = \emptyset$. This implies that $V(t) = V(x(t)) > R$ for all $t \in [0, t_{\max})$. Inequality (2.1) implies that $\dot{V}(t) \leq 0$ for almost all $t \in [0, t_{\max})$, which implies $V(x(t)) \leq V(x_0)$ for all $t \in [0, t_{\max})$. Since $V \in C^2(\Re^n; \Re_+)$ is radially unbounded, it follows from the inequality $V(x(t)) \leq V(x_0)$ that the solution $x(t)$ is bounded on $[0, t_{\max})$. Thus, standard theory of ordinary differential equations implies that $t_{\max} = +\infty$. Let $\delta(x_0) > 0$ be defined by $\delta(x_0) := \min\{W(x) : R \leq V(x) \leq V(x_0)\}$. Indeed, notice that positivity of $\delta(x_0) > 0$ is a direct consequence of continuity of $W$ and compactness of the set $\{x \in \Re^n : R \leq V(x) \leq V(x_0)\}$. Differential inequality (2.1) in conjunction with $V(x(t)) \leq V(x_0)$ and $V(t) = V(x(t)) > R$ implies that $\dot{V}(t) \leq -\delta(x_0)$ for almost all $t \in [0, +\infty)$. Consequently, we obtain $R < V(x(t)) \leq V(x_0) - \delta(x_0) t$, for all $t \geq 0$, which is a contradiction.

Since $A := \{t \in [0, t_{max}) : V(x(t)) \leq R\} \neq \emptyset$, we define $t_1 = \inf A$. Continuity of $V$ implies that $t_1 > 0$ and $V(x(t_1)) = R$. Moreover, positive invariance of the compact set $S = \{x \in \Re^n : V(x) \leq R\}$ implies that the solution exists for all $t \geq 0$ and satisfies $V(x(t)) \leq R$ for all $t \geq t_1$. Using an argument similar to the one used for the case $A = \emptyset$ we are in a position to establish that $V(x(t)) \leq V(x_0)$ for all $t \in [0, t_1)$ and $t_1 \leq \frac{V(x_0) - R}{\delta(x_0)}$, where $\delta(x_0) := \min\{W(x) : R \leq V(x) \leq V(x_0)\}$.

Finally, define $T(x_0) := 0$ for all $x_0 \in S$ and $T(x_0) := \frac{V(x_0) - R}{\delta(x_0)}$ for all $x_0 \notin S$, where $\delta(x_0) := \min\{W(x) : R \leq V(x) \leq V(x_0)\}$. The above analysis guarantees that $V(x(t)) \leq R$ for all $t \geq T(x_0)$ and that $V(x(t)) \leq \max(V(x_0), R)$ for all $t \geq 0$. Continuity of the function $T : \Re^n \to \Re_+$ is a direct consequence of continuity of $V, W$ and the fact that the level sets of $V$ are compact sets. The proof is complete. ◁

**Proof of the claim that system (3.16) is forward complete:** First notice that for every initial condition and every input the component $x(t)$ of the solution $(x(t), z(t))$ of system (3.16) is defined for all $t \geq 0$. This follows from the fact that (1.1), (1.2) is forward complete.

By virtue of the definition $\dot{W}(z, x, u) = \nabla W(z) \tilde{k}(z, h(x), u) - \lambda(z, h(x), u) \nabla W(z) Q(z) (\nabla W(z))'$ and the fact that $W$ is a radially unbounded function, we guarantee the existence of a continuous function $\tilde{c} : H(A) \times U \to [1, +\infty)$ such that:

$$\dot{W}(z, x, u) \leq \tilde{c}(h(x), u), \text{ for all } u \in U, (z, x) \in \Re^n \times \Re^n \text{ with } W(z) \leq a + 1 \tag{A.1}$$

Inequality (A.1) in conjunction with (3.15) shows that there exists a continuous function $\hat{c} : H(A) \times U \to [1, +\infty)$ such that:

$$\dot{W}(z, x, u) \leq \hat{c}(h(x), u) W(z), \text{ for all } u \in U, (z, x) \in \Re^n \times \Re^n \tag{A.2}$$

The differential inequality (A.2) shows that the solution $(x(t), z(t))$ of system (3.16) satisfies the following inequality for almost all $t \geq 0$ for which the solution exists:

$$\dot{W}(t) \leq \beta(t) W(t) \tag{A.3}$$

where $W(t) = W(z(t))$ and $\beta(t) := \hat{c}(h(x(t)), u(t))$. The differential inequality (A.3) shows that $W(z(t)) \leq \exp\left(\int_0^t \beta(s) ds\right) W(z(0))$ for all $t \geq 0$ for which the solution exists. Moreover, the inequality $W(z(t)) \leq \exp\left(\int_0^t \beta(s) ds\right) W(z(0))$ shows that $W(z(t))$ remains bounded on bounded intervals of time. Using the fact that $W$ is a radially unbounded function and a standard contradiction argument, we conclude that the solution $(x(t), z(t))$ of system (3.16) is defined for all $t \geq 0$. The proof is complete. ◁